\documentclass[11pt]{amsart}
\usepackage[utf8]{inputenc}
\usepackage[T1]{fontenc}
\usepackage{lmodern}
\usepackage{amsmath,amssymb,amsthm}
\usepackage{booktabs}
\usepackage{xcolor}
\usepackage[colorlinks=true,linkcolor=blue!50!black,citecolor=blue!50!black,urlcolor=blue!50!black]{hyperref}
\usepackage[margin=1.15in]{geometry}

\numberwithin{equation}{section}
\newtheorem{theorem}{Theorem}[section]
\newtheorem{lemma}[theorem]{Lemma}
\newtheorem{proposition}[theorem]{Proposition}
\newtheorem{corollary}[theorem]{Corollary}
\newtheorem*{theoremW}{Theorem W \textup{(Ward)}}
\newtheorem*{theoremZ}{Theorem Z \textup{(Zsigmondy)}}
\newtheorem*{corollary*}{Corollary}
\theoremstyle{definition}
\newtheorem{definition}[theorem]{Definition}
\newtheorem{rulex}[theorem]{Rule}
\theoremstyle{remark}
\newtheorem{remark}[theorem]{Remark}

\newcommand{\I}{I}
\newcommand{\sups}{\operatorname{sup}}
\newcommand{\ord}{\operatorname{ord}}
\newcommand{\pmin}{p_{\min}}
\newcommand{\psup}{p_{\sup}}
\newcommand{\Bmin}{B_{\min}}
\newcommand{\Tsafe}{T_{\mathrm{safe}}}

\title{Every friend of 10 has at least eleven distinct prime factors}
\author{Filip Fifka}
\email{filip.fifka@gmail.com}

\date{14 November 2026}

\begin{document}

\maketitle

\begin{abstract}
A \emph{friend of $10$} is a positive integer $N\neq 10$ such that $\sigma(N)/N=\sigma(10)/10=9/5$.
Whether one exists is an open problem. Thackeray proved that a friend of $10$ has at least $10$ distinct prime
factors. We describe a computer-assisted proof that it has at least $11$. The proof is an
exhaustive search over sets of primes, carried out by a recursive certifier written in PARI/GP,
which extends the open-source project \texttt{caioross/solitary10}. The new ingredients are a split
of ``almost tight'' nodes that avoids enumerating the smallest unknown prime, a treatment of nodes
with three or more unknown primes by a saturation dichotomy, a counting rule for unfactored
cofactors of cyclotomic values, and the correction of a soundness gap in the earlier closure rule.
The certification run for $\omega=10$ visited $157{,}210{,}962$ nodes in $11.4$ hours on $12$
threads. As a consequence, a friend of $10$ would be larger than $10^{38}$. Every prime the
computation relies on is proved, the largest ones by ECPP certificates distributed with the code,
which is available at \url{http://elliptic.group/files/solitary10-gp.tar.bz2}.
\end{abstract}

\section{Introduction}

For a positive integer $n$ let $\sigma(n)$ be the sum of its divisors and
$\I(n)=\sigma(n)/n$ its \emph{abundancy index}. Two distinct integers are \emph{friends} if they
have the same abundancy index; an integer with no friend is \emph{solitary}. The smallest integer
whose status is unknown is $10$, with $\I(10)=9/5$.

\subsection*{Earlier results}
Ward~\cite{ward} showed that a friend $N$ of $10$ is an odd square, that its least prime factor is
$5$, that $3\nmid N$ and that $N$ has at least $6$ distinct prime factors. Thackeray~\cite{thackeray}
proved that $N$ has at least $10$ distinct prime factors, using a bound on the
exponents of the known primes (his Corollary~6 and Proposition~9, which we will extend in this paper). Necessary
conditions and bounds on the individual prime divisors of $N$ are studied in
\cite{mandal1,mandal3,mandal4}; \cite{mandal2} treats the analogous question for $20$. To our
knowledge (as of 11 November 2026) no bound beyond Thackeray's has been announced in the literature.

The project \texttt{solitary10}~\cite{repo} is an open-source Python implementation of a recursive
certifier for the statement ``no friend of $10$ has exactly $k$ distinct prime factors''. It
re-derives Ward's theorems, eliminates $k\le 8$, and documents the node at which the method stops
for $k=9$. Our computational work starts from that code and its documentation, after correcting a closure rule
that was unsound (see Section~\ref{sec:f5}).

\subsection*{Results}
Write $\omega(N)$ for the number of distinct prime factors of $N$.

\begin{theorem}\label{thm:main}
If $\sigma(N)/N=9/5$ and $N\ne 10$, then $\omega(N)\ge 11$.
\end{theorem}

\noindent Combined with an exhaustive bounded search (Section~\ref{sec:size}) this gives:

\begin{corollary*}[see Corollary~\ref{cor:size}]
Every friend of $10$ is larger than
$10^{38}$.
\end{corollary*}

The proof of Theorem~\ref{thm:main} consists of the application of mathematical rules described in
Sections~\ref{sec:setting}--\ref{sec:f5}, each of which removes only configurations that contain no
friend of $10$, and of a computation (Section~\ref{sec:comp}) that applies them until nothing is
left for $\omega(N)=k$, $1\le k\le 10$.

Section~\ref{sec:notation} fixes the notation and collects the elementary facts used throughout,
with proofs. Section~\ref{sec:setting} defines the nodes of the search tree, what it means for a
rule to be sound, and proves the index lemma, which is the source of every size bound.
Section~\ref{sec:old} recalls the rules inherited from \cite{repo}. Sections~\ref{sec:b7}
and~\ref{sec:b8} state and prove the new rules. Section~\ref{sec:f5} describes the corrected
closure rule. Sections~\ref{sec:comp}--\ref{sec:verif} describe the computations and their
verification.

\section{Notation and elementary facts}\label{sec:notation}

\subsection{Notation}
Throughout, $p,q,r,u$ and $U$ (with or without indices) denote primes, and $\ell$ is used for the
lower bound of a node as defined in Section~\ref{sec:setting}. When $\ell$ denotes a prime (in the
budgets of Section~\ref{sec:old}) we say so explicitly.

\medskip
\begin{center}
\begin{tabular}{ll}
\toprule
Symbol & Meaning\\
\midrule
$\sigma(n)$ & sum of the positive divisors of $n$\\
$\I(n)$ & abundancy index $\sigma(n)/n$\\
$\omega(n)$ & number of distinct prime factors of $n$\\
$v_r(x)$ & exponent of the prime $r$ in the rational number $x\neq0$\\
$\ord_q(p)$ & multiplicative order of $p$ modulo the prime $q\nmid p$\\
$\Phi_d(x)$ & $d$-th cyclotomic polynomial\\
$\sups(p)$ & $p/(p-1)$, the limit of $\I(p^a)$ as $a\to\infty$\\
$N$ & a friend of $10$: $\sigma(N)/N=9/5$, $N\ne10$\\
$S$ & the set of prime factors of $N$; $N=\prod_{p\in S}p^{a_p}$\\
$a_p$ & the exponent of $p\in S$ in $N$\\
$\lfloor x\rfloor$ & the largest integer $\le x$\\
$(x,y]$ & the set of integers $n$ with $x<n\le y$ (or the primes in it, when listing primes)\\
\bottomrule
\end{tabular}
\end{center}
\medskip

\noindent For a prime $p$ and $a\ge1$ we write $\I(p^a)=\sigma(p^a)/p^a$. Note that $\sups(p)$ is
defined for every integer $p\ge2$ by the same formula, and we use it for integers as well.

\subsection{Facts about the abundancy index}

\begin{lemma}\label{lem:I}
Let $n,m\ge1$ and let $p\ge2$ be an integer.
\begin{enumerate}
\item[(i)] $\I(n)=\sum_{d\mid n}1/d$.
\item[(ii)] If $\gcd(n,m)=1$ then $\I(nm)=\I(n)\I(m)$.
\item[(iii)] For $a\ge1$, $\I(p^a)=1+p^{-1}+\dots+p^{-a}=\dfrac{p^{a+1}-1}{p^a(p-1)}$; it is
strictly increasing in $a$ and strictly smaller than $\sups(p)$.
\item[(iv)] For a real number $B>1$, let $\sups(B)=1+\dfrac1{B-1}$ and, by a harmless abuse of
notation, $\I(B^2)=1+\dfrac1B+\dfrac1{B^2}$; when $B$ is a prime these agree with the definitions
above, and they are only ever evaluated at bounds for primes. Both are strictly decreasing
functions of $B$, and $\I(B^2)<\sups(B)$.
\end{enumerate}
\end{lemma}
\begin{proof}
(i) The map $d\mapsto n/d$ is a bijection of the set of divisors of $n$ onto itself, so
$\sigma(n)=\sum_{d\mid n}d=\sum_{d\mid n}n/d=n\sum_{d\mid n}1/d$.
(ii) When $\gcd(n,m)=1$, every divisor of $nm$ is uniquely a product $de$ with $d\mid n$,
$e\mid m$, so $\sum_{f\mid nm}1/f=\sum_{d\mid n}1/d\cdot\sum_{e\mid m}1/e$.
(iii) The divisors of $p^a$ are $1,p,\dots,p^a$; the closed form is the geometric sum. Each
increase of $a$ adds the positive term $p^{-a-1}$, and the infinite sum is $1/(1-1/p)=\sups(p)$.
(iv) Both expressions decrease as $B$ grows, and
$\sups(B)=\sum_{j\ge0}B^{-j}>1+1/B+1/B^2$.
\end{proof}

\subsection{Ward's theorem and its consequences}

\begin{theoremW}[\cite{ward}]
Let $N$ be a friend of $10$. Then $N$ is an odd square, its least prime factor is $5$, and
$3\nmid N$. Equivalently, $N=\prod_{p\in S}p^{a_p}$ with $\min(S)=5$ and every $a_p$ even and at least $2$.
\end{theoremW}

\noindent We use Theorem~W as an established result. The following consequences are elementary.

\begin{lemma}[Product equation, valuations, closure]\label{lem:val}
Let $N$ be a friend of $10$. Then
\begin{equation}\label{eq:prod}
\frac95=\prod_{p\in S}\I\bigl(p^{a_p}\bigr),
\end{equation}
\begin{equation}\label{eq:val}
v_3\bigl(\sigma(N)\bigr)=2,\qquad v_5\bigl(\sigma(N)\bigr)=a_5-1,\qquad
v_r\bigl(\sigma(N)\bigr)=a_r\quad(r\in S,\ r\neq5),
\end{equation}
and every prime factor of $\sigma(p^{a_p})$, for every $p\in S$, lies in $S\cup\{3\}$.
\end{lemma}
\begin{proof}
Equation \eqref{eq:prod} is $\I(N)=9/5$ together with Lemma~\ref{lem:I}(ii). From
$5\sigma(N)=9N$ we get $\sigma(N)=3^2\cdot5^{-1}\cdot N$, and since $3\notin S$ this gives
\eqref{eq:val}. Every prime factor of $\sigma(N)=9N/5$ is therefore $3$ or a prime of $N$. Finally
$\sigma$ is multiplicative, so $\sigma(N)=\prod_{p\in S}\sigma(p^{a_p})$, and each prime factor of
a factor $\sigma(p^{a_p})$ divides $\sigma(N)$.
\end{proof}

\subsection{Cyclotomic values}

The polynomial identity $x^n-1=\prod_{d\mid n}\Phi_d(x)$ gives, for a prime $p$ and $a\ge1$,
\begin{equation}\label{eq:cyc}
\sigma(p^a)=\frac{p^{a+1}-1}{p-1}=\prod_{d\mid a+1,\ d>1}\Phi_d(p).
\end{equation}
For a friend of $10$, every $a_p$ is even, so $a_p+1$ is odd and only odd $d$ occur in
\eqref{eq:cyc}. A prime $q$ dividing $\Phi_d(p)$ is called \emph{primitive} if $\ord_q(p)=d$.

\begin{lemma}\label{lem:cycdiv}
Let $p$ be a prime, $d\ge1$, and $q\ne p$ a prime dividing $\Phi_d(p)$. Then $d=\ord_q(p)\cdot q^e$
for some $e\ge0$. In particular: if $q$ is primitive then $\equiv1\pmod d$; if $q$ is not
primitive then $q\mid d$; and a prime dividing both $\Phi_d(p)$ and $\Phi_{d'}(p)$ with $d\ne d'$
divides $dd'$.
\end{lemma}
\begin{proof}
The first statement is a standard property of cyclotomic polynomials (see e.g.~\cite{zsigmondy}).
If $\ord_q(p)=d$ then $d\mid q-1$ by Fermat's little theorem. If $e\ge1$ then $q\mid d$. If $q$
divides $\Phi_d(p)$ and $\Phi_{d'}(p)$ with $d\ne d'$, then $d=\ord_q(p)q^e$ and $d'=\ord_q(p)q^{e'}$ with $e\ne e'$. 
Since $e,e' \geq 0 $ are integers, it follows that $\max({e,e'})\geq 1 $
so $q$ divides at least the larger of $d,d'$.
\end{proof}

\begin{theoremZ}[\cite{zsigmondy}]
Let $a>b\ge1$ be coprime integers and $n\ge1$. Then $a^n-b^n$ has a prime factor that divides no
$a^k-b^k$ with $1\le k<n$ (a \emph{primitive prime divisor}), except in the following cases:
$n=1$ and $a-b=1$; $n=2$ and $a+b$ is a power of $2$; $(a,b,n)=(2,1,6)$.
\end{theoremZ}

\noindent We use it in the following form.

\begin{lemma}\label{lem:zsig}
Let $p\ge5$ be a prime and $d\ge3$ an integer. Then $\Phi_d(p)$ has a primitive prime factor,
that is, a prime $q$ with $q\mid\Phi_d(p)$ and $\ord_q(p)=d$.
\end{lemma}
\begin{proof}
Apply Theorem~Z with $a=p$, $b=1$ and $n=d$; these are coprime and $a>b$. None of the exceptional
cases occurs: $n=d\ge3$ excludes the first two, and $a=p\ge5$ excludes $(2,1,6)$. So $p^d-1$ has a
prime factor $q$ that divides no $p^k-1$ with $1\le k<d$. In particular $q\ne p$, since $p$ does not
divide $p^d-1$, so $\ord_q(p)$ is defined.

 From $q\mid p^d-1$ we get $\ord_q(p)\mid d$. If $\ord_q(p)=k<d$, then
$q\mid p^k-1$, which is excluded. Hence $\ord_q(p)=d$.

By the identity $x^n-1=\prod_{k\mid n}\Phi_k(x)$ used in
\eqref{eq:cyc}, $p^d-1=\prod_{k\mid d}\Phi_k(p)$, so the prime $q$ divides $\Phi_k(p)$ for some
$k\mid d$. If $k<d$, then $\Phi_k(p)$ divides $p^k-1$, so $q\mid p^k-1$, which is again excluded.
Hence $k=d$, and $q$ is a primitive prime factor of $\Phi_d(p)$.
\end{proof}

By Lemma~\ref{lem:cycdiv}, primitive primes for different $d$ are different. Hence:

\begin{corollary}\label{cor:zsig}
For a friend $N$ of $10$ and $p\in S$, each divisor $d>1$ of $a_p+1$ provides a prime $q_d \in S$ such that
$q_d  \equiv1\pmod d$, and different divisors provide different primes.
\end{corollary}
\begin{proof}
Lemma~\ref{lem:zsig} gives a primitive prime $q_d$ of $\Phi_d(p)$; it divides $\sigma(p^{a_p})$ by
\eqref{eq:cyc}, so $q_d\in S\cup\{3\}$ by Lemma~\ref{lem:val}; $q_d\equiv1\pmod d$ with $d\ge3$
excludes $q_d=3$; and $\ord_{q_d}(p)=d$ determines $d$.
\end{proof}

\subsection{A lemma on fractions}
The following trivial fact is used several times; we state it to make the later proofs explicit.

\begin{lemma}[Denominator lemma]\label{lem:den}
Let $x=A/D$ with integers $A$ and $D\ge1$ (not necessarily in lowest terms). If $x>1$ then
$x\ge1+1/D$.
\end{lemma}
\begin{proof}
$x>1$ means $A>D$, so $A\ge D+1$.
\end{proof}

\section{Nodes, soundness and the index lemma}\label{sec:setting}

\subsection{Nodes}
The certifier proves, for one value of $k$ at a time, that no friend of $10$ has exactly $k$
distinct prime factors. It does so by exploring a tree of \emph{nodes}, each describing a set of
possible friends.

\begin{definition}[Node]\label{def:node}
A \emph{node} $\mathcal N=(C,s,\ell,F,m,M,\kappa)$ consists of:
\begin{itemize}
\item a finite set $C$ of primes (with at least $5\in C$), the \emph{known} primes;
\item an integer $s\ge0$, the number of \emph{unknown} primes;
\item an integer $\ell\ge5$, the \emph{lower bound};
\item a partial map $F$ on $C$ assigning \emph{fixed} exponents $F(p)$;
\item for the primes $p\in C$ without a fixed exponent (the \emph{free} primes), a lower bound
$m(p)\ge2$ (default $2$) and possibly an upper bound $M(p)$;
\item a \emph{constraint} $\kappa$, which is either empty or one of the constraints defined in
Sections~\ref{sec:b7} and~\ref{sec:b8}.
\end{itemize}
A friend $N$ of $10$ \emph{belongs to} $\mathcal N$ if: $C\subseteq S$ and $|S|=|C|+s$; the
exponents satisfy $a_p=F(p)$ for fixed $p$ and $m(p)\le a_p\le M(p)$ for free $p$; the constraint
$\kappa$ holds; and the \textbf{invariant} holds: \emph{every prime of $S$ that is at most $\ell$
lies in $C$}.
\end{definition}

\noindent For a friend belonging to $\mathcal N$ we write the primes of $S\setminus C$ in
increasing order as $U_1<U_2<\dots<U_s$ and their exponents as $a_1,\dots,a_s$. By the invariant,
$U_1>\ell$. The \emph{root node} for a given $k$ is $(\{5\},k-1,5,\emptyset,\dots)$: every friend
with $\omega(N)=k$ belongs to it, since $5$ is its least prime (Theorem~W).

\subsection{The bounds attached to a node}
For a node, let $\underline a_p$ be $F(p)$ if $p$ is fixed and $m(p)$ otherwise, and let
\[
\overline\I_p=\begin{cases}\I\bigl(p^{F(p)}\bigr)&\text{if $p$ is fixed},\\
\I\bigl(p^{M(p)}\bigr)&\text{if $p$ is free with an upper bound},\\
\sups(p)&\text{otherwise.}\end{cases}
\]
Define
\begin{equation}\label{eq:pminsup}
\pmin=\prod_{p\in C}\I\bigl(p^{\underline a_p}\bigr),\qquad \psup=\prod_{p\in C}\overline\I_p,\qquad
R=\frac{9/5}{\psup},\qquad \delta=R-1 .
\end{equation}
For a friend $N$ belonging to the node, write $\I_C=\prod_{p\in C}\I(p^{a_p})$ for the true
contribution of the known primes.

\begin{lemma}\label{lem:bracket}
For a friend $N$ belonging to a node,
\begin{equation}\label{eq:E}
\pmin\le\I_C\le\psup\qquad\text{and}\qquad\frac95=\I_C\prod_{i=1}^s\I\bigl(U_i^{a_i}\bigr).
\end{equation}
Moreover $R\le 45/31<2$.
\end{lemma}
\begin{proof}
By Lemma~\ref{lem:I}(iii), $\I(p^{a})$ is increasing in $a$ and below $\sups(p)$, so
$\I(p^{\underline a_p})\le\I(p^{a_p})\le\overline\I_p$ for each $p\in C$; multiply. The equation
is \eqref{eq:prod} split into the known and the unknown primes. Since $5\in C$ with exponent at
least $2$, $\psup\ge\I(5^2)=31/25$ (the other factors are greater or equal to $1$), so $R\le\frac95\cdot\frac{25}{31}=\frac{45}{31}$.
\end{proof}

\noindent The number $\delta$ measures the gap the unknown primes must fill. A node is called
\emph{almost tight} when $\delta$ is very small.

\subsection{Rules and soundness}

\begin{definition}[Soundness rule]
A \emph{rule} replaces a node $\mathcal N$ by finitely many nodes $\mathcal N_1,\dots,\mathcal N_t$
(possibly $t=0$), its \emph{children}. The rule is \emph{sound} if every friend of $10$ that
belongs to $\mathcal N$ belongs to at least one child.
\end{definition}

\begin{proposition}\label{prop:tree}
Suppose the certifier, starting from the root node for $k$, applies only sound rules, and every
node it creates is eventually replaced (in particular every node with $s=0$ is decided by an exact
search, see Section~\ref{sec:old}). If the process ends with no node left, then no friend of $10$
has exactly $k$ distinct prime factors.
\end{proposition}
\begin{proof}
Suppose $N$ is a friend with $\omega(N)=k$. It belongs to the root. By soundness, whenever a node
to which $N$ belongs is replaced, $N$ belongs to one of its children. So at every stage some node
in the current set contains $N$, and the set is never empty.
\end{proof}

\noindent If some node cannot be treated by any rule, the program records a \emph{refusal}; a run
with a refusal proves nothing. A rule that kills a node (replaces it by no children) is sound
exactly when no friend belongs to that node.

\subsection{The index lemma}

\begin{lemma}[Index lemma]\label{lem:index}
\textup{(a)} For a friend belonging to a node,
\[
\prod_{i=1}^s\sups(U_i)\ >\ \frac{9/5}{\I_C}\ \ge\ R .
\]
\textup{(b)} Let $j\ge1$ and let $2\le V_1\le V_2\le\dots\le V_j$ be integers with
$\prod_{i=1}^j\sups(V_i)>\rho$, where $1<\rho\le2$. Then
\[
V_1\ <\ 1+\frac{2j}{\rho-1}.
\]
\end{lemma}
\begin{proof}
(a) By \eqref{eq:E}, $\prod_i\I(U_i^{a_i})=\frac95/\I_C$, and $\I(U_i^{a_i})<\sups(U_i)$ by
Lemma~\ref{lem:I}(iii); the second inequality is $\I_C\le\psup$.

(b) Since $\sups$ is decreasing (Lemma~\ref{lem:I}(iv)) and $V_1$ is the smallest, each factor is
at most $\sups(V_1)$, so $\sups(V_1)^j\ge\prod_i\sups(V_i)>\rho$. Put $x=1/(V_1-1)>0$, so that
$\sups(V_1)=1+x$. Using $1+x\le e^x$ we get $e^{jx}\ge(1+x)^j>\rho$, hence $jx>\ln\rho$. For
$\rho>1$ we have $\ln\rho\ge1-1/\rho=(\rho-1)/\rho$ (the inequality $\ln y\ge1-1/y$ holds for all
$y>0$), and $\rho\le2$ gives $(\rho-1)/\rho\ge(\rho-1)/2$. So $jx>(\rho-1)/2$, that is,
$V_1-1=1/x<2j/(\rho-1)$.
\end{proof}

\begin{remark}
By Lemma~\ref{lem:bracket}, part (b) always applies with $\rho=R$ and $j=s$, giving
$U_1<1+2s/\delta$. This bound is finite, but very large when $\delta$ is small: that is the
difficulty Blocks~7 and~8 address.
\end{remark}

\section{The rules inherited from \texttt{solitary10}}\label{sec:old}

We recall the rules of \cite{repo} (its ``Blocks 3--6''). They are restated and proved in the
review document distributed with the code; here we state what each one does and the fact it relies
on. Below, $\ell'$ denotes a prime (not the lower bound of a node).

\begin{enumerate}
\item[(R1)] \emph{Lower index bound.} If $s\ge1$ and $\pmin\ge9/5$, the node is empty, since every
unknown contributes a factor $>1$ in \eqref{eq:E}.
\item[(R2)] \emph{Complete sets.} If $s=0$, the finitely many exponent vectors allowed by the
closure rule (R3) are tested exactly against \eqref{eq:prod}.
\item[(R3)] \emph{Closure.} For a fixed exponent $F(p)$, the value $\sigma(p^{F(p)})$ is computed
piece by piece by \eqref{eq:cyc} and factored (up to a size limit, $160$ bits in our runs). By
Lemma~\ref{lem:val}, each prime factor outside $C\cup\{3\}$ is an unknown prime of $S$; it is
\emph{pinned}, that is, moved into $C$, and $s$ decreases by one. A pin $u\le\ell$ contradicts the
invariant, so the node is empty. A composite part that cannot be factored within the limit is an
\emph{unfactored cofactor} (Section~\ref{sec:cof}).
\item[(R4)] \emph{Budgets of $3$ and $5$.} For $\ell'\in\{3,5\}$ the number $\ell'-1$ has no odd
divisor above $1$, while $\ord_{\ell'}(q)$ must divide the odd number $a_q+1$ whenever
$\ell'\mid\sigma(q^{a_q})$; so only primes $q\equiv1\pmod{\ell'}$ contribute, and for them
$v_{\ell'}(\sigma(q^{a}))=v_{\ell'}(a+1)$ (lifting the exponent). By \eqref{eq:val},
\[
\sum_{q\in S,\ q\equiv1\,(3)}v_3(a_q+1)=2,\qquad \sum_{q\in S,\ q\equiv1\,(5)}v_5(a_q+1)=a_5-1,
\]
and each power $\ell'^j\mid a_q+1$ requires its own primitive prime of $\Phi_{\ell'^j}(q)$ in $S$
(Corollary~\ref{cor:zsig}).
\item[(R5)] \emph{Last unknown.} If $s=1$, some prime $\ell'$ divides $a_U+1$, so $\Phi_{\ell'}(U)$
divides $\sigma(N)$; by Lemma~\ref{lem:cycdiv} its prime factors are $\ell'$ (to the power at most
$1$) and primes $\equiv1\pmod{\ell'}$, which must be known. So
$\Phi_{\ell'}(U)=\ell'^{\varepsilon}\prod_{q\in C,\,q\equiv1\,(\ell')}q^{e_q}$, and since
$\Phi_{\ell'}$ is increasing, each right-hand side gives at most one $U$. Once $U$ or the $e_q$ are
bounded there are finitely many candidates.
\item[(R6)] \emph{Exponent bound} (Thackeray~\cite{thackeray}, Corollary~6 and Proposition~9). See
Proposition~\ref{prop:cor6} below.
\item[(R7)] \emph{Exponent branching.} A free exponent $a_p$ (even, greater or equal to $m(p)$) is replaced by the
children $a_p=m(p),m(p)+2,\dots,A-2$ (fixed) and one \emph{tail} child with $m(p):=A$. This is a
partition of the possible values, hence sound for every choice of $A$.
\item[(R8)] \emph{Listing the smallest unknown.} If nothing else applies, the certifier creates,
for each prime $q>\ell$, $q\notin C$, in increasing order, the child $(C\cup\{q\},s-1,q,\dots)$
in which $q=U_1$; this is sound because no prime of $N$ lies strictly between $\ell$ and $U_1$.
The listing stops at the first $q$ with
$\psup\,\sups(q)\prod(\text{the $s-1$ next sups})\le9/5$; this quantity decreases with $q$, so by
Lemma~\ref{lem:index}(a) no later $q$ can be $U_1$ either.
\end{enumerate}

For the next statement we need some notation. Fix a prime $r\in C$ (so $r\ge5$ is odd). Put
$P=C\cup\{3\}$. For $j\ge1$ let $W_j$ be the number of primes $t\in P$, $t\ne r$, with
$t\equiv1\pmod{r^j}$, and for $q\in C$ let $W_j^{(q)}$ be the same number with $q$ left out. For a
sequence $W=(W_j)$ and an integer $c\ge0$ let
\[
\lambda(W,c)=\max\bigl\{v\ge0:\ W_j+c\ \ge\ v-j+1\ \text{ for all }1\le j\le v\bigr\}.
\]
For an integer $x$ prime to $r$ let $w_r(x)=v_r(x^{r-1}-1)\ge1$. For $t\ge1$ let
\[
m_t(r)=\min_{2\le y\le r-1}\bigl(y^{\,r^{t-1}}\bmod r^t\bigr),\qquad
\omega_r(L)=\max\{t\ge1:\ m_t(r)\le L\}.
\]
Note that $m_1(r)=2$, so $\omega_r(L)\ge1$ for every $L\ge2$.
 
\begin{proposition}[Exponent bound, after Thackeray]\label{prop:cor6}
Let $r\in C$ be a free prime and let $L_1,\dots,L_s$ be integers. Let $c_5=1$ and
$c_r=0$ for $r\ne5$, and
\begin{align*}
A\;=\;c_r\;&+\sum_{q\in C\ \mathrm{fixed}}v_r\bigl(\sigma(q^{F(q)})\bigr)
\;+\sum_{\substack{q\in C\ \mathrm{free}\\ q\ne r}}\beta(q)\\
&+\;\sum_{i=1}^{s}\bigl(\omega_r(L_i)+\lambda(W,s-1)\bigr),
\end{align*}
where, with $\lambda_q=\min\bigl(\lambda(\cdot,s),\ \lfloor\log_r(M(q)+1)\rfloor\bigr)$ (the second
term only if $q$ has an upper bound $M(q)$),
\[
\beta(q)=\begin{cases}
\lambda_q\ \text{ computed with } W^{(q)} & \text{if } q\equiv1\pmod r,\\
w_r(q)+\lambda_q\ \text{ computed with } W & \text{if } q\not\equiv1\pmod r \text{ and } \ord_r(q) \text{ is odd},\\
0 & \text{if } \ord_r(q) \text{ is even}.
\end{cases}
\]
Then every friend belonging to the node with $U_i\le L_i$ $(1\le i\le s)$ has $a_r\le A$. (The
bound holds friend by friend; it is not required that every friend of the node satisfies
$U_i\le L_i$.)
\end{proposition}
\begin{proof}
Let $N$ be such a friend.
 
\emph{Step 0: the valuation identity.} By \eqref{eq:val}, $v_r(\sigma(N))=a_r-c_r$. Since
$\sigma(r^{a_r})\equiv1\pmod r$ and $\sigma$ is multiplicative,
\begin{equation}\label{eq:arsum}
a_r=c_r+\sum_{q\in S,\ q\ne r}v_r\bigl(\sigma(q^{a_q})\bigr).
\end{equation}
For a known prime with a fixed exponent the term is computed exactly. We need to bound the term
of a free known prime $q$ by $\beta(q)$, and the term of each unknown $U_i$ by
$\omega_r(L_i)+\lambda(W,s-1)$.
 
\emph{Step 1: the valuation of one term.} Let $q\ne r$ be a prime and $a\ge2$ even, and let
$e=\ord_r(q)$. Since $\sigma(q^a)=(q^{a+1}-1)/(q-1)$ and $r$ is odd, lifting the exponent gives:
\begin{itemize}
\item if $q\equiv1\pmod r$: $v_r(q^{a+1}-1)=v_r(q-1)+v_r(a+1)$, so $v_r(\sigma(q^a))=v_r(a+1)$;
\item if $e>1$ and $e\mid a+1$: $v_r(q-1)=0$ and
$v_r(q^{a+1}-1)=v_r(q^e-1)+v_r\bigl((a+1)/e\bigr)$. Here $v_r(q^e-1)=w_r(q)$ and
$v_r((a+1)/e)=v_r(a+1)$, because $e$ divides $r-1$, which is prime to $r$. So
$v_r(\sigma(q^a))=w_r(q)+v_r(a+1)$;
\item if $e>1$ and $e\nmid a+1$: $r\nmid q^{a+1}-1$ and $v_r(\sigma(q^a))=0$.
\end{itemize}
The exponent $a$ is even (Theorem~W), so $a+1$ is odd, and an even $e$ never divides it: a prime
$q$ of even order modulo $r$ contributes $0$. In all cases
\begin{equation}\label{eq:oneterm}
v_r\bigl(\sigma(q^a)\bigr)\ \le\ \varepsilon(q)\,w_r(q)+v_r(a+1),\qquad
\varepsilon(q)=\begin{cases}1&\text{if }q\not\equiv1\!\!\pmod r\text{ and }e\text{ is odd},\\0&\text{otherwise,}\end{cases}
\end{equation}
and the left-hand side is $0$ when $e$ is even.
 
\emph{Step 2: the levels $v_r(a+1)$.} Let $q\in S$, $q\ne r$, and $j=v_r(a_q+1)$. For
$1\le i\le j$ the number $r^i$ divides $a_q+1$, so $\Phi_{r^i}(q)$ divides $\sigma(q^{a_q})$ by
\eqref{eq:cyc}. By Lemma~\ref{lem:zsig}, $\Phi_{r^i}(q)$ has a prime factor $t_i$ with $\ord_{t_i}(q)=r^i$;
hence $t_i\equiv1\pmod{r^i}$, the primes $t_1,\dots,t_j$ are distinct, none is $q$ or $r$, and
each lies in $S\cup\{3\}$ (Lemma~\ref{lem:val}). Fix $i\le j$. The $j-i+1$ primes $t_i,\dots,t_j$
are all $\equiv1\pmod{r^i}$. Those that are known lie in $P\setminus\{q,r\}$, so there are at most
$W_i$ of them ($W_i^{(q)}$ if $q\in C$). The others are unknown primes different from $q$: at
most $s$ of them if $q$ is known, at most $s-1$ if $q$ is itself unknown. Therefore
$j-i+1\le W_i+c$ for every $i\le j$, with $c=s$ or $s-1$, that is $j\le\lambda(W,c)$. Moreover
$r^j\le a_q+1\le M(q)+1$ when $q$ has an upper bound. This gives $v_r(a_q+1)\le\lambda_q$ for a
free known prime and $v_r(a_i+1)\le\lambda(W,s-1)$ for an unknown, and with \eqref{eq:oneterm} the
bound $\beta(q)$ for a free known prime.
 
\emph{Step 3: the term $w_r(U)$ of an unknown.} Let $U$ be an unknown with $U\le L$. If
$\varepsilon(U)=0$ there is nothing to prove, since $\omega_r(L)\ge1$. Otherwise
$U\not\equiv1\pmod r$; let $t=w_r(U)$, so $U^{r-1}\equiv1\pmod{r^t}$. The congruence
$x^{r-1}\equiv1\pmod{r^t}$ has exactly $r-1$ solutions modulo $r^t$, one in each non-zero class
modulo $r$, namely $y^{r^{t-1}}\bmod r^t$ for $1\le y\le r-1$. Since $U\not\equiv1\pmod r$, the
residue of $U$ modulo $r^t$ is one of those with $y\ge2$, so $U\ge m_t(r)$, hence $m_t(r)\le L$.
The sequence $m_t(r)$ is non-decreasing in $t$, because a solution modulo $r^{t+1}$ reduces to a
solution modulo $r^t$ in the same class modulo $r$, and the residue can only grow. So
$t\le\omega_r(L)$.
 
Putting Steps 1--3 into \eqref{eq:arsum} gives $a_r\le A$.
\end{proof}
 
\noindent This is essentially Corollary~6 and Proposition~9 of \cite{thackeray}, with one change: Step~3 is
applied to each unknown separately, so each unknown $U_i$ may have its own bound $L_i$, where
\cite{thackeray} uses one bound for all. 
\noindent We use the method of \cite[Corollary~6]{thackeray}, but not its original statement. The primes counted
in Step~2 are the primitive prime factors $t_i$ of $\Phi_{r^i}(q)$;
they divide $\sigma(N)$, and Step~2 places them in $S\cup\{3\}$ by Lemma~\ref{lem:val} without any
assumption on their size. The bounds $L_i$ are used only in Step~3, for the unknown prime factors
of $N$. (In \cite[Corollary~6]{thackeray} the hypothesis is stated for the prime factors $q>r$ of
$n$; it should be stated for every prime $q>r$, or for the prime factors $q>r$ of $\sigma(n)$, as
H.~R.~Thackeray has confirmed. For a friend of $10$ and $r\ge5$ the three
readings coincide, since the prime factors of $\sigma(N)=9N/5$ are those of $N$ together with
$3<r$.)
\noindent  Besides Theorem~Z and Lemma~\ref{lem:val} (Steps~0 and~2), the proof uses only that the exponents
are even (Step~1) and that $r$ is odd. Of the node it uses only $C$, $s$, the fixed exponents and
the upper bounds $M(q)$; it does not use the lower bound $\ell$, the constraint $\kappa$ or any
index bound, so the bound holds at every node and for every choice of the $L_i$.
\noindent The bound holds for every $r$ but the program uses it only for 
$r \leq 10^6$, because computing $m_t(r)$ costs time proportional to $r$. For larger $r$ the rule is not applied.
The bound is also tested against true values during the computation (Section~\ref{sec:verif}).

\section{Block 7: almost tight nodes}\label{sec:b7}

Throughout this section, $\mathcal N$ is a node with $s\ge1$ and $\psup<9/5$ at which the listing
(R8) would otherwise be used.

\subsection{Two thresholds}
\begin{definition}
\[
B'=\max\Bigl\{B\ge2:\ \psup\,\I(B^2)>\tfrac95\Bigr\},\qquad
\Bmin=\max\Bigl\{U\ge2:\ \psup\,\sups(U)\ge\tfrac95\Bigr\}.
\]
\end{definition}
\noindent Both sets are defined through the functions $B\mapsto\I(B^2)=1+1/B+1/B^2$ and
$U\mapsto\sups(U)=1+1/(U-1)$, which are strictly decreasing for real $B,U>1$ and tend to $1$
(Lemma~\ref{lem:I}(iv)). Hence the left-hand sides $\psup\,\I(B^2)$ and $\psup\,\sups(U)$ decrease
to $\psup<9/5$, so both sets are bounded. They are non-empty, because $5\in C$ gives
$\psup\ge\I(5^2)=31/25$, so that $\psup\,\I(2^2)=\tfrac74\psup>\tfrac95$ and
$\psup\,\sups(2)=2\psup>\tfrac95$. Finally, since the left-hand sides decrease, the defining
inequality holds for every $2\le B\le B'$ and every $2\le U\le\Bmin$.

\begin{lemma}\label{lem:bmin}
$\Bmin=\lfloor1+1/\delta\rfloor$, and $B'\le\Bmin$.
\end{lemma}
\begin{proof}
$\psup\sups(U)\ge9/5$ means $\sups(U)\ge R$, that is $1+\frac1{U-1}\ge1+\delta$, that is
$U\le1+1/\delta$. For the second claim, $\I(B^2)<\sups(B)$ (Lemma~\ref{lem:I}(iv)), so
$\psup\I(B'^2)>9/5$ implies $\psup\sups(B')>9/5$, hence $B'\le\Bmin$.
\end{proof}

\noindent Since $\I(B^2)=1+1/B+1/B^2$ and $\sups(B)=1+1/B+1/B^2+1/B^3+\cdots$ differ only by
terms of order $B^{-3}$, the gap $(B',\Bmin]$ is tiny; in the census of $k=9$ it was empty at
every node.

\subsection{The split}
\begin{definition}[Constraint ``hi'']
The constraint $\kappa=h$, for an integer $h$, means: \emph{the smallest unknown prime satisfies
$U_1\le h$}.
\end{definition}

\begin{rulex}[Split]\label{rule:split}
Replace $\mathcal N$ by the following children:
\begin{enumerate}
\item[(a)] $\mathcal N$ with the added constraint $\kappa=B'$;
\item[(b)] for each prime $q\in(B',\Bmin]$, $q\notin C$, the child of (R8) in which $U_1=q$;
\item[(c)] $\mathcal N$ with $\ell$ replaced by $\Bmin$.
\end{enumerate}
\end{rulex}
\begin{proof}[Soundness]
Let $N$ belong to $\mathcal N$. Its smallest unknown $U_1$ lies in exactly one of $(\ell,B']$,
$(B',\Bmin]$ and $(\Bmin,\infty)$. In the first case $N$ belongs to (a). In the second it belongs
to the child (b) with $q=U_1$. In the third case every unknown exceeds $\Bmin$, so every prime of
$S$ that is smaller or equal to $\Bmin$ is known: the invariant holds with $\ell=\Bmin$ and $N$ belongs to (c). The
argument uses no property of $B'$ and $\Bmin$ other than $\ell\le B'\le\Bmin$, so floating-point
errors in computing them could only affect efficiency; the program nevertheless checks both
thresholds exactly.
\end{proof}

\noindent\emph{A node that already carries a constraint.} If $\mathcal N$ already carries
$\kappa=h$ (alone, or as part of a saturation constraint, Definition~\ref{def:sat}), every child
of the split keeps it: (b) is listed only over $(B',\min(h,\Bmin)]$, child (c) is omitted when
$h\le\Bmin$, and the split is not applied at all when $h\le B'$. Intersecting each child with
$U_1\le h$ still covers every friend of $\mathcal N$, since every such friend satisfies $U_1\le h$.

\subsection{The hi rule}
The point of child (a) is that the constraint $U_1\le B'$ gives a lower bound on the factor of
$U_1$, which kills the child without listing $U_1$.

\begin{lemma}[Hi rule]\label{lem:hi}
Let $\mathcal N$ carry the constraint $\kappa=h$.
\begin{enumerate}
\item[(i)] Every friend belonging to $\mathcal N$ satisfies $\pmin\,\I(h^2)\le9/5$, with strict
inequality if $s\ge2$. Hence $\mathcal N$ is empty if $\pmin\I(h^2)>9/5$, or if
$\pmin\I(h^2)=9/5$ and $s\ge2$, or if $h\le\ell$.
\item[(ii)] Suppose a rule pins a set $P$ of unknown primes. If no prime of $P$ is smaller or equal to $ h$, the
child may keep the constraint $\kappa=h$. If moreover $P$ contains all $s$ unknowns, the child is
empty.
\end{enumerate}
\end{lemma}
\begin{proof}
(i) Let $N$ belong to $\mathcal N$. Since $a_1\ge2$ and $U_1\le h$, Lemma~\ref{lem:I}(iii) and (iv)
give $\I(U_1^{a_1})\ge\I(U_1^2)\ge\I(h^2)$. By \eqref{eq:E},
\[
\frac95=\I_C\,\I(U_1^{a_1})\prod_{i\ge2}\I(U_i^{a_i})\ \ge\ \pmin\,\I(h^2)\prod_{i\ge2}\I(U_i^{a_i}),
\]
and the last product is $1$ if $s=1$ and $>1$ if $s\ge2$. If $h\le\ell$ then $U_1\le h\le\ell$
contradicts the invariant.

(ii) The smallest unknown $U_1$ satisfies $U_1 \le h$. If no pinned prime is smaller or equal to $h$, then $U_1\notin P$, so
$U_1$ is still unknown in the child and is the smallest of the remaining unknowns; the constraint
holds. If $P$ contains all the unknowns, then $U_1\in P$ and $U_1\le h$, a contradiction.
\end{proof}

\noindent In child (a), $h=B'$ and, by definition of $B'$, $\psup\,\I(B'^2)>9/5$. The certifier
treats $V=\I(B'^2)$ as an extra factor: it branches the free exponents (R7) and, in each tail,
$\pmin$ grows. Because $\psup V>9/5$, there is room for $\pmin V$ to exceed $9/5$, at which point
Lemma~\ref{lem:hi}(i) kills the node; exact children are treated by the closure. The range
$(\ell,B']$, which contains almost all of the millions of candidates for $U_1$, is thus never
listed.

\subsection{Branch (c) with two unknowns}

\begin{lemma}\label{lem:c2}
Let $s=2$ and suppose every unknown exceeds $\Bmin$. Write $R=n/m$ in lowest terms. Then
\[
U_1\ <\ 1+\frac4\delta=:L_1\qquad\text{and}\qquad U_2\ <\ 1+m\,U_1\ \le\ 1+m\,L_1 .
\]
\end{lemma}
\begin{proof}
The bound on $U_1$ is Lemma~\ref{lem:index}(b) with $j=2$, $\rho=R$, using part (a). Since
$U_1>\Bmin$, the definition of $\Bmin$ gives $\psup\sups(U_1)<9/5$, that is $\sups(U_1)<R$. Put
$\rho_2=R/\sups(U_1)>1$. By Lemma~\ref{lem:index}(a), $\sups(U_1)\sups(U_2)>R$, that is
$\sups(U_2)>\rho_2$. Now
\[
\rho_2=\frac{n}{m}\cdot\frac{U_1-1}{U_1}=\frac{n(U_1-1)}{mU_1}
\]
is a fraction $>1$ with denominator $mU_1$, so $\rho_2-1\ge1/(mU_1)$ by Lemma~\ref{lem:den}.
Finally $\sups(U_2)>\rho_2$ means $1/(U_2-1)>\rho_2-1$, that is
$U_2<1+1/(\rho_2-1)\le1+mU_1$.
\end{proof}

\noindent The program uses the weaker bound $U_2<1+2m\,L_1$. So in branch (c) with $s=2$ every
unknown is bounded, and Proposition~\ref{prop:cor6} bounds the exponent of every free prime
$r\le10^6$ of $C$.

\subsection{Narrowing}
\begin{lemma}[Narrowing]\label{lem:narrow}
Let $s=2$, and suppose the closure (R3) of the fixed exponents produced no pin and at least one
unfactored cofactor $c$. Then both unknown primes divide $c$; in particular $U_1\mid c$.
\end{lemma}
\begin{proof}
By Lemma~\ref{lem:cofac} below, $c$ has at least two distinct prime factors, all of them unknown
primes of $S$. There are only two unknown primes, so they are exactly those factors.
\end{proof}

\noindent The listing (R8) then runs only over the prime factors of $c$ below the stop, found by
trial division, instead of over every prime.

\subsection{The safe chain and the band}
\begin{lemma}[Safe chain]\label{lem:chain}
Let $s\ge3$ and define $\eta$ by
\[
1+\eta=\frac{9/5}{\psup\,\sups(\ell+1)^{\,s-1}} .
\]
If $\eta>0$, then every unknown prime of every friend belonging to the node is $<1+1/\eta$.
\end{lemma}
\begin{proof}
All unknowns exceed $\ell$, so by Lemma~\ref{lem:I}(iv) each of $U_1,\dots,U_{s-1}$ has
$\sups(U_i)\le\sups(\ell+1)$. Lemma~\ref{lem:index}(a) gives
$\sups(U_s)\prod_{i<s}\sups(U_i)>R$, hence
$\sups(U_s)>R/\sups(\ell+1)^{s-1}=1+\eta$, that is $U_s<1+1/\eta$. The other unknowns are smaller
than $U_s$.
\end{proof}

\noindent Let $\Tsafe$ be the least $\ell$ for which $\eta>0$; it is about $(s-1)/\delta$. Branch
(c) is cut into two sound children: ``$U_1\le\Tsafe$'' (constraint $\kappa=\Tsafe$) and ``every
unknown exceeds $\Tsafe$'' ($\ell:=\Tsafe$). In the second, Lemma~\ref{lem:chain} bounds every
unknown and Proposition~\ref{prop:cor6} applies. What remains to be listed is the \emph{band}
\[
U_1\in(\Bmin,\Tsafe]\qquad(s\ge3),
\]
where the index alone does not bound $U_2,\dots,U_s$. At $k=10$ a band holds between $10^4$ and
$10^8$ primes.

\section{Block 8: the band}\label{sec:b8}

Throughout this section, $\mathcal N$ is a node with $s\ge2$ unknowns (from
Section~\ref{sec:dich} on, $s\ge3$) and $\psup<9/5$. The only exception is the lower lift
(Rule~\ref{rule:lift}), which the program applies at every node with $s\ge1$; its proof uses
neither restriction.

\subsection{The lower lift}
\begin{rulex}[Lower lift]\label{rule:lift}
Let $B''$ be the largest integer $B$ with $\pmin\,\I(B^2)\ge9/5$ if $s\ge2$ (resp.
$\pmin\,\I(B^2)>9/5$ if $s=1$), or $B''=\ell$ if there is none. Replace $\ell$ by $\max(\ell,B'')$.
\end{rulex}
\begin{proof}[Soundness]
Let $N$ belong to $\mathcal N$. As in the proof of Lemma~\ref{lem:hi}(i),
$9/5\ge\pmin\,\I(U_1^2)\prod_{i\ge2}\I(U_i^{a_i})$, so $\pmin\I(U_1^2)<9/5$ if $s\ge2$ and
$\pmin\I(U_1^2)\le9/5$ if $s=1$. If some unknown prime were $\le B''$, then $U_1\le B''$, and since
$\I(B^2)$ decreases, $\pmin\I(U_1^2)\ge\pmin\I(B''^2)$, which is at least $9/5$ (resp. strictly greater than $9/5$): a
contradiction. So every unknown exceeds $B''$, every prime of $S$ that is at most $\max(\ell,B'')$ is
known, and the invariant holds with the new value.
\end{proof}

\subsection{Deepening}
Let $T=\lfloor1+2s/\delta\rfloor$, so that $U_1\le T$ by Lemma~\ref{lem:index}(b) (if the node
carries $\kappa=h$ with $h<T$, the program uses $T=h$ instead, which is also an upper bound for
$U_1$; everything below holds with either value), let $n_f$ be
the number of free primes in $C$, and $\tau=1/(n_fT^3)$. Before the region step, the certifier
applies exponent branching (R7) to the free prime with the largest ratio
$\overline\I_p/\I(p^{\underline a_p})$ until every free prime satisfies
\[
\frac{\overline\I_p}{\I\bigl(p^{\underline a_p}\bigr)}-1\ \le\ \tau .
\]
Each step is an instance of (R7), so this is sound whatever $\tau$ is. Its effect is that $\pmin$
and $\psup$ agree to within a relative error of order $T^{-3}$, which is what makes the test of
Section~\ref{sec:sat} effective.

\subsection{The saturation dichotomy}\label{sec:dich}
For each $p\in C$ choose a number $x_p\ge\overline\I_p$; the program takes
$x_p=\sups(p)$ when $\sups(p)/\overline\I_p-1\le\tau$, and $x_p=\overline\I_p$ otherwise (so
$x_p=\sups(p)$ for a free prime without an upper bound, and $x_p=\I(p^{M(p)})$ for a free prime
with an upper bound $M(p)$ far from $\sups(p)$). Put
\[
p^*=\prod_{p\in C}x_p\ \ge\ \psup,\qquad R^*=\frac{9/5}{p^*}=\frac{n^*}{m^*}\ \ (\text{lowest terms}).
\]
The purpose of $p^*$ is to keep the denominator $m^*$ small; replacing $\overline\I_p$ by the
simple fraction $p/(p-1)$ removes the huge denominators of exact values $\I(p^a)$ with large $a$.
If $R^*\le1$ the region step below is not used. Otherwise $1<R^*\le R<2$.

\begin{definition}[The two regions]
A friend belonging to $\mathcal N$ is \emph{saturated} if
\[
\prod_{i=1}^{s-1}\sups(U_i)\ \ge\ R^*,
\]
and \emph{unsaturated} otherwise. In words: the $s-1$ smallest unknown primes, at their largest
possible contributions, already reach $R^*$, so the largest unknown may be arbitrarily large.
\end{definition}

\begin{lemma}[Bounds in the unsaturated region]\label{lem:unsat}
Define $L_1=T$ and, for $1\le t\le s-1$,
\[
L_{t+1}=1+2(s-t)\,m^*\prod_{i=1}^{t}L_i .
\]
Every unsaturated friend belonging to $\mathcal N$ satisfies $U_i\le L_i$ for $1\le i\le s$.
\end{lemma}
\begin{proof}
Fix an unsaturated friend. First, $U_1\le T=L_1$ by Lemma~\ref{lem:index}(b). We show, for
$1\le t\le s-1$,
\begin{equation}\label{eq:unsatstep}
U_{t+1}\ <\ 1+2(s-t)\,m^*\prod_{i=1}^{t}U_i ;
\end{equation}
since the right-hand side increases with $U_1,\dots,U_t$, induction on $t$ then gives $U_i\le L_i$.

\emph{Step 1: a ratio larger than one.} Put
\[
R^*_t=\frac{R^*}{\prod_{i=1}^t\sups(U_i)} .
\]
Each factor $\sups(U_i)$ exceeds $1$, so
$\prod_{i\le t}\sups(U_i)\le\prod_{i\le s-1}\sups(U_i)<R^*$ (the friend is unsaturated). Hence
$R^*_t>1$.

\emph{Step 2: how much larger.} Since $\sups(U)=U/(U-1)$,
\[
R^*_t=\frac{n^*}{m^*}\cdot\prod_{i=1}^t\frac{U_i-1}{U_i}=\frac{n^*\prod_{i\le t}(U_i-1)}{m^*\prod_{i\le t}U_i},
\]
a fraction larger than $1$ with denominator $D=m^*\prod_{i\le t}U_i$. By Lemma~\ref{lem:den},
$R^*_t-1\ge1/D$.

\emph{Step 3: the remaining unknowns.} By Lemma~\ref{lem:index}(a) and $\I_C\le\psup\le p^*$,
\[
\prod_{i=1}^s\sups(U_i)\ >\ \frac{9/5}{\I_C}\ \ge\ \frac{9/5}{p^*}=R^*,\qquad\text{so}\qquad
\prod_{i=t+1}^{s}\sups(U_i)\ >\ R^*_t .
\]
The $s-t$ integers $U_{t+1}<\dots<U_s$ and $\rho=R^*_t$, which satisfies $1<\rho\le R^*<2$, meet the
hypotheses of Lemma~\ref{lem:index}(b) with $j=s-t$. Therefore
\[
U_{t+1}\ <\ 1+\frac{2(s-t)}{R^*_t-1}\ \le\ 1+2(s-t)\,D ,
\]
which is \eqref{eq:unsatstep}.
\end{proof}

\subsection{The region step}
\begin{definition}[Saturation constraint]\label{def:sat}
For an integer $h$, an integer $k\ge2$ and a rational $\rho_0>1$, the constraint
$\kappa=[h,k,\rho_0]$ means:
\[
U_1\le h\qquad\text{and}\qquad\prod_{i=1}^{k}\sups(U_i)\ \ge\ \rho_0 .
\]
It is used only when $k=s-1$. Whenever a rule pins an unknown prime (so that the set of unknowns
changes), the child drops the second part of the constraint, which is about the unknowns of the
parent; the first part, $U_1\le h$, is kept as long as no pinned prime is smaller or equal to $h$
(Lemma~\ref{lem:hi}(ii)). Rules that pin nothing (exponent branching, deepening, the split of
Rule~\ref{rule:split}, narrowing) pass the constraint on unchanged to every child; the split and
narrowing list their children by Rule~\ref{rule:sat} below. Dropping a constraint, or a part of
it, enlarges the set of friends belonging to a node, so it is always sound.
\end{definition}

\begin{rulex}[Region step]\label{rule:region}
Let $r\in C$ be a free prime with $r\le10^6$, and let $A$ be the bound of
Proposition~\ref{prop:cor6} computed with the bounds $L_1,\dots,L_s$ of Lemma~\ref{lem:unsat}.
Replace $\mathcal N$ by:
\begin{itemize}
\item the children with $a_r$ fixed to each admissible even value in $[m(r),A]$;
\item one tail child with $m(r)$ raised to the least admissible value strictly above$A$, carrying the
constraint $[h,s-1,R^*]$, where
\[
h=\min\Bigl(T,\ \max\bigl\{U\ge2:\ \sups(U)^{s-1}>R^*\bigr\}\Bigr).
\]
\end{itemize}
\end{rulex}
\begin{proof}[Soundness]
The children partition the possible values of $a_r$, so it suffices to show that every friend with
$a_r>A$ satisfies the constraint. Such a friend cannot be unsaturated: if it were, then
$U_i\le L_i$ for all $i$ by Lemma~\ref{lem:unsat}, and Proposition~\ref{prop:cor6} would give
$a_r\le A$. So it is saturated: $\prod_{i<s}\sups(U_i)\ge R^*$. Moreover, $s-1\ge2$ and
$U_1<U_2$, so $\sups(U_1)>\sups(U_2)$ and
\[
\sups(U_1)^{s-1}\ >\ \prod_{i=1}^{s-1}\sups(U_i)\ \ge\ R^*,
\]
which, $\sups$ being decreasing, means $U_1\le\max\{U:\sups(U)^{s-1}>R^*\}$. Also $U_1\le T$. So
$U_1\le h$. The number $R^*$ is a fixed rational, so the constraint stays true in every
descendant with the same unknown primes.
\end{proof}

\subsection{Listing under saturation}\label{sec:sat}
Let $\mathcal N$ carry the constraint $[h,k,\rho_0]$ with $k=s-1\ge2$. Its smallest unknown is
listed, $U_1=q$ for the primes $q\in(\ell,h]$, $q\notin C$. For each $q$ put
\[
\rho=\rho(q)=\frac{\rho_0}{\sups(q)}.
\]
Dividing the constraint by $\sups(U_1)=\sups(q)$ gives, for every friend in the child $U_1=q$,
\begin{equation}\label{eq:rho}
\prod_{i=2}^{k}\sups(U_i)\ \ge\ \rho .
\end{equation}
A prime $q$ with $\pmin\,\I(q^2)\ge9/5$ is skipped, because its child is empty: by
Lemma~\ref{lem:bracket}, $9/5=\I_C\,\I(q^{a_1})\prod_{i\ge2}\I(U_i^{a_i})>\pmin\,\I(q^2)$, since
$a_1\ge2$ and $s\ge2$. If $\rho\le1$ the child gets no constraint (nothing is claimed). Otherwise:

\begin{rulex}[Saturated listing]\label{rule:sat}
Assume $\rho>1$.
\begin{enumerate}
\item[(i)] If $k\ge3$, let $H=\max\{U\ge2:\sups(U)^{k-1}>\rho\}$. If $H\le q$ the child is
empty; otherwise it carries the constraint $[H,k-1,\rho]$.
\item[(ii)] If $k=2$ (that is, $s=3$), let
\[
H(q)=\Bigl\lfloor\frac{\rho}{\rho-1}\Bigr\rfloor .
\]
If $H(q)\le q$, or if
\begin{equation}\label{eq:kill}
\pmin\,\I(q^2)\,\I\bigl(H(q)^2\bigr)\ \ge\ \frac95 ,
\end{equation}
the child is empty. Otherwise it carries the constraint $\kappa=H(q)$ (Lemma~\ref{lem:hi}).
\end{enumerate}
\end{rulex}
\begin{proof}[Soundness]
(i) In the child, the constrained unknowns are $U_2<\dots<U_k$, which are $k-1\ge2$ integers. As
in the proof of Rule~\ref{rule:region}, $\sups(U_2)^{k-1}>\prod_{i=2}^k\sups(U_i)\ge\rho$, so
$U_2\le H$. Since $U_2>U_1=q$, the case $H\le q$ is impossible. The new constraint
$[H,k-1,\rho]$ is exactly \eqref{eq:rho} together with $U_2\le H$, where $U_2$ is now the
smallest unknown and $k-1=(s-1)-1$.

(ii) Here \eqref{eq:rho} reads $\sups(U_2)\ge\rho$. Since $\rho>1$,
\[
\frac{U_2}{U_2-1}\ge\rho\iff U_2\ge\rho U_2-\rho\iff U_2(\rho-1)\le\rho\iff U_2\le\frac{\rho}{\rho-1},
\]
so $U_2\le H(q)$, and $H(q)\le q$ is impossible as before. Now use \eqref{eq:E} with the three
unknowns $q=U_1$, $U_2$, $U_3$ and bound each factor from below: $\I_C\ge\pmin$,
$\I(q^{a_1})\ge\I(q^2)$, $\I(U_2^{a_2})\ge\I(U_2^2)\ge\I(H(q)^2)$ (Lemma~\ref{lem:I}), and
$\I(U_3^{a_3})>1$. This gives
\[
\frac95\ =\ \I_C\,\I(q^{a_1})\,\I(U_2^{a_2})\,\I(U_3^{a_3})\ >\ \pmin\,\I(q^2)\,\I\bigl(H(q)^2\bigr).
\]
So if \eqref{eq:kill} holds, no friend belongs to the child. The remaining constraint $U_2\le H(q)$
is a hi constraint for the child, whose smallest unknown is $U_2$.
\end{proof}

\begin{lemma}[Integer form of the test]\label{lem:int}
Write $\rho_0=n/m$ and $\pmin/(9/5)=a/b$ with positive integers, and put $d=n(q-1)-mq$. Then
$\rho>1$ if and only if $d>0$; in that case
\[
H(q)=\Bigl\lfloor\frac{n(q-1)}{d}\Bigr\rfloor,
\]
and \eqref{eq:kill} holds if and only if
\[
a\,(q^2+q+1)\,(H^2+H+1)\ \ge\ b\,q^2H^2\qquad(H=H(q)).
\]
\end{lemma}
\begin{proof}
$\rho=\frac nm\cdot\frac{q-1}{q}=\frac{n(q-1)}{mq}$, so $\rho-1=\frac{d}{mq}$, which is positive
if and only if $d>0$, and then $\rho/(\rho-1)=n(q-1)/d$. For the test, $\I(x^2)=(x^2+x+1)/x^2$,
so \eqref{eq:kill} reads $\frac ab\cdot\frac{(q^2+q+1)(H^2+H+1)}{q^2H^2}\ge1$; multiply by the
positive number $bq^2H^2$.
\end{proof}

\begin{remark}[Why almost every $q$ dies]\label{rem:why}
To first order, with exponents $2$, the constraint asks $1/q+1/U_2\gtrsim\delta^*$ (with
$R^*=1+\delta^*$), while \eqref{eq:E} asks $1/q+1/U_2+1/U_3\lesssim\delta$, and after deepening
$\delta$ and $\delta^*$ agree to within $O(T^{-3})$. So $U_2$ is squeezed into a window whose width
is of the order of the third-order terms; \eqref{eq:kill} checks whether that window still contains
an integer, and it almost never does. By Lemma~\ref{lem:int} each test costs a few integer
multiplications. In the certification run for $k=10$, $5.16\times10^9$ candidates $q$ were
discarded in this way.
\end{remark}

\subsection{Counting unfactored cofactors}\label{sec:cof}

\begin{definition}[Unfactored cofactor]
When the closure (R3) factors a piece $\Phi_d(p)$ of $\sigma(p^{F(p)})$, it removes the primes
below $10^5$, the primes of $C\cup\{3\}$ and the primes pinned so far. If the remaining part $c>1$
is composite, is
not a prime power, and is too large to factor within the budget, $c$ is called an
\emph{unfactored cofactor}.
\end{definition}

\begin{lemma}\label{lem:cofac}
An unfactored cofactor $c$ has at least two distinct prime factors, and each of them is an unknown
prime of $S$ (an element of $S\setminus C$).
\end{lemma}
\begin{proof}
A composite integer that is not a prime power has at least two distinct prime factors. Each prime
factor of $c$ divides $\sigma(p^{F(p)})$, so it lies in $S\cup\{3\}$ (Lemma~\ref{lem:val}); it is
not $3$ and not in $C$, since those primes were removed from $c$.
\end{proof}

\begin{rulex}[Coprime cofactors]\label{rule:coprime}
Let $c_1,\dots,c_t$ be the unfactored cofactors produced at a node, let $P$ be the product of all
primes pinned at that node, and let $\mathcal F\subseteq\{c_1,\dots,c_t\}$ be a family of cofactors
that are pairwise coprime and each coprime to $P$. If
\[
\#\{\text{pins}\}+2\,\#\mathcal F\ >\ s,
\]
the node is empty.
\end{rulex}
\begin{proof}[Soundness]
For a friend belonging to the node, the pins are distinct unknown primes. By
Lemma~\ref{lem:cofac}, each $c\in\mathcal F$ contributes at least two distinct unknown primes;
they are not pins ($\gcd(c,P)=1$), and two different members of $\mathcal F$ share none (they are
coprime). So $N$ has at least $\#\{\text{pins}\}+2\#\mathcal F$ distinct unknown primes, but it has
exactly $s$.
\end{proof}

\noindent The program applies the rule to a whole node (step~1 of the partition) and also, in the
budget cases (R4), to the pins and cofactors of each single assignment of the budget. The proof is
the same in both cases: only the pins and the cofactors of the set under consideration are used.

\begin{lemma}[Tight family]\label{lem:tight}
If no prime was pinned and $2\,\#\mathcal F=s$, then every unknown prime divides
$\prod_{c\in\mathcal F}c$, and each $c\in\mathcal F$ has exactly two distinct prime factors.
\end{lemma}
\begin{proof}
By the proof of Rule~\ref{rule:coprime}, the members of $\mathcal F$ contain at least $2\#\mathcal
F=s$ distinct unknown primes. There are exactly $s$ unknowns, so these are all of them, and no
member can have a third distinct prime factor.
\end{proof}

\noindent Lemma~\ref{lem:tight} extends narrowing (Lemma~\ref{lem:narrow}, which is the case
$s=2$, $\#\mathcal F=1$) to $s=4,6,\dots$: the smallest unknown is searched only among the prime
factors of $\prod_{c\in\mathcal F}c$ below the index stop.


\section{A correction to the closure rule}\label{sec:f5}

The closure of \cite{repo}, and our port of it, discarded a node when $\#\{\text{pins}\}+e>s$,
where $e$ is twice the largest number of unfactored cofactors coming from a single known prime. A
cofactor is stripped of the primes pinned \emph{before} it is computed, but a prime pinned
\emph{later}, from the closure of another fixed exponent, may divide it. Its two ``further''
unknown primes then include that pin, and the count is too large.

\begin{proposition}\label{prop:f5}
The rule above is unsound. Let $p_3=1091$, so that $u=\sigma(p_3^2)=1191373$ is prime, and let
$q_1=50038757$, so that $\sigma(q_1^2)=u\cdot2101673659$ with $2101673659$ prime. With a factoring
budget below $52$ bits and the fixed exponents processed in the order $q_1,p_3$, the node with
$C=\{5,p_3,q_1\}$, $F(q_1)=F(p_3)=2$ and $s=2$ is discarded, although the two unknown primes
needed, $u$ and $2101673659$, fit into $s=2$.
\end{proposition}
\begin{proof}
Processing $q_1$ first, $\sigma(q_1^2)$ (of $52$ bits) cannot be factored within the budget and
becomes an unfactored cofactor: $e=2$. Processing $p_3$ next pins $u$. The rule counts
$1+2=3>2$. But the cofactor $u\cdot2101673659$ is divisible by the pin $u$, so only one further
unknown prime is needed.
\end{proof}

\noindent A planted target with this structure (Section~\ref{sec:verif}) was indeed lost by the
earlier program. The corrected rule is Rule~\ref{rule:coprime} applied to the cofactors of each
single prime, with $P$ the product of \emph{all} pins of the node. Its premise holds: the cofactors
of one prime $p$ come from different pieces $\Phi_d(p)$ and $\Phi_{d'}(p)$, which by
Lemma~\ref{lem:cycdiv} share only primes dividing $dd'$; these are below $10^5$ (in our runs
$d<10^5$ is enforced) and were removed, so the cofactors are pairwise coprime. The same correction
is made in the budget cases (R4). With the real budget of $160$ bits the faulty situation needs a
coincidence between two unrelated cyclotomic values. The node counts for $k\le8$ are identical
before and after the correction, and for $k=9$ the corrected program again produces exactly
$36{,}554{,}260$ nodes, the count of the original run; so the gap did not affect any of the runs
for $k\le9$.

\section{The computation}\label{sec:comp}

The long runs were made with PARI/GP~2.17.2 on a $12$-thread machine with ECC memory; the tests of
Section~\ref{sec:verif} and the smaller timings of Section~\ref{sec:size} with PARI/GP~2.15.4 on a
$2$-thread machine (the searches for $B=10^{37}$ and $10^{38}$ ran on the $12$-thread machine). The
program is \texttt{core/omega\_k.gp} with the driver \texttt{experiments/certify\_omega\_k.gp}.
The top of the tree is cut into tasks that are run in parallel; a task exceeding a time limit is abandoned and expanded into
smaller tasks, and every node is counted once. A run is accepted only if no task is pending or
failed, no node is refused and no friend is found.

\begin{table}[h]
\centering
\begin{tabular}{rrl}
\toprule
$k$ & nodes & remarks\\
\midrule
$1$--$7$ & $1,1,3,5,7,55,459$ & \\
$8$ & $26{,}222$ & \\
$9$ & $36{,}554{,}260$ & $1$ h $28$ min on $12$ threads (Block 7 only)\\
$10$ & $157{,}210{,}962$ & $11$ h $23$ min on $12$ threads, $15$ rounds (Blocks 7 and 8)\\
\bottomrule
\end{tabular}
\caption{Size of the certified trees.}
\end{table}

For $k=10$ the settings are: cofactors up to $160$ bits are factored; $\Phi_d(p)$ is evaluated for
$\varphi(d)\le2000$ and at most $4000$ bits; the exponent bound is used for known primes below
$10^6$; Block~8 is applied to bands of more than $1000$ primes. Before the certification run, a
\emph{census} run with the same rules (which records, instead of listing, any node whose listing
exceeds $10^4$ primes) finished in $8$ h $49$ min on $12$ threads with $149{,}823{,}612$ nodes, no
refusal and four recorded nodes, all with $C=\{5,7,11,13,29,19531\}$, $a_5=6$, $s=4$ and a band of
$16{,}690$ primes. These four nodes were then run separately and closed in $45$ to $106$ minutes
each ($1.0$ to $4.1$ million nodes, $8.1\times10^8$ discarded primes each). The certification run
itself (same settings, with the primality policy of Section~\ref{sec:primes}) finished after $15$
rounds and $40{,}988$ s with $157{,}210{,}962$ nodes, no pending or failed task, no refused node and
no friend; $5.16\times10^9$ primes were discarded by the arithmetic test of Rule~\ref{rule:sat} and
$8.1\times10^6$ listings were narrowed to the prime factors of cofactors. Rounds $9$--$14$ each split
the same few tasks again; they are pieces of the four heavy nodes above.

\subsection{Primality}\label{sec:primes}
A composite cofactor wrongly declared prime would create a false pin, so primality matters in one
direction. Factorisations and multiplicative orders use proven primes. A cofactor of at most $400$
bits is proved prime; a larger one is accepted if it passes the Baillie--PSW test and $20$ rounds of
Miller--Rabin with random bases, and is written to a list, to be proved separately. For $k=9$ the
list has $436$ primes, and all of them are proved: $432$ by the APRCL test (\texttt{isprime}), and
the four largest (of about $5230$, $6684$, $7283$ and $17512$ bits; two of them are
$\Phi_{125}(q)/1255$ and $\Phi_{125}(q')/5$), which are beyond the practical range of APRCL, by the
ECPP implementation of PARI/GP (\texttt{primecert}). The four ECPP certificates were checked with
\texttt{primecertisvalid} and are distributed with the code, so they can be re-checked independently.
For $k=10$ the list has $4{,}328$ primes of $401$ to $3{,}596$ bits, and all $4{,}328$ are proved by
the APRCL test. Thus the runs for $k\le10$ rest on no unproved primality claim.

\subsection{Implementation notes}
(i) In the ``last unknown'' rule (R5), a prime $\ell'\mid r-1$ of tens of millions makes
$\Phi_{\ell'}(U_{\max})$ a number of nearly $10^9$ bits. For such $\ell'$ the solutions $U$ are found
by scanning the primes $U$ of the (short) admissible interval and computing $v_q(\Phi_{\ell'}(U))$
modulo powers of $q$, with the size comparison $\Phi_{\ell'}(U)>U^{\ell'-1}$. (ii) The time limit of
a task is tested between nodes and, inside a saturated listing, every $8192$ listed primes. The
partial results of an abandoned task are discarded and the task is expanded again from its nodes,
so where the limit is tested does not affect soundness. (iii) Narrowing needs every prime factor of
a cofactor $c$ below the index stop $D$. These are found with PARI's \texttt{factor(c,D)}, which
returns all prime factors of $c$ below $D$ (we checked this against trial division, also for
$D>2^{32}$).

\section{A lower bound for the least friend}\label{sec:size}

\subsection{The bounded search}
The project \cite{repo} contains an exhaustive bounded search. For a bound $B$ and a number $k$ of
distinct primes, it enumerates the sets $S$ of $k$ primes with $\min S=5$ and
$\prod_{p\in S}p^2\le B$, and for each set the exponent vectors $(a_p)$ (even, $a_p \ge2$) with
$\prod_{p\in S}p^{a_p}\le B$ (the cyclotomic closure is not used in this search); it prunes with exact index bounds as in
Section~\ref{sec:setting}, and tests every remaining vector, called a \emph{signature}, against
$\sigma(N)/N=9/5$. Since every exponent is at least $2$, a friend $N$ satisfies
$N\ge\bigl(\prod_{p\in S}p\bigr)^2$, so only finitely many $k$ can occur below $B$. Starting that
search at $k=10$ (Thackeray's bound), \cite{repo} reports that the least friend of $10$ exceeds
$10^{32}$.

Theorem~\ref{thm:main} alone gives only $N\ge(5\cdot7\cdots41)^2\approx2.6\times10^{27}$. But the
case $k=10$ is by far the largest part of the bounded search, and without it the search reaches
much further.

\begin{corollary}\label{cor:size}
A friend of $10$ is larger than $10^{38}$.
\end{corollary}
\begin{proof}
Let $N$ be a friend of $10$ and $k=\omega(N)$. By Theorem~\ref{thm:main}, $k\ge11$.

\emph{Large $k$.} The product of the $15$ smallest primes greater or equal to $5$ is approximately $1.03\times10^{41}>10^{38}$. Since every exponent of $N$ is at least $2$
and its primes are at least the $k$ smallest primes greater or equal to $5$, a friend with $k\ge15$ is larger than
$10^{38}$.

\emph{$11\le k\le14$.} The bounded search for $B=10^{38}$ finds no friend for $k=11,12,13,14$.
It visits $9{,}632{,}870{,}504$, $17{,}633{,}755$, $338$ and $5$ nodes and tests
$9{,}843{,}884{,}172$, $10{,}534{,}449$, $0$ and $0$ signatures respectively,
$9{,}854{,}418{,}621$ signatures in all, in $36{,}520$ s on $12$ threads ($36{,}439$ s of it at
$k=11$). (For $k=14$ the search is needed: $(5\cdot7\cdots53)^2\approx2.95\times10^{37}<10^{38}$.)
\end{proof}

\begin{table}[h]
\centering
\begin{tabular}{rrr}
\toprule
$B$ & signatures tested, $k\ge 11$ & time\\
\midrule
$10^{32}$ & $3{,}633$ & $0.2$ s$^{}$\\
$10^{33}$ & $139{,}674$ & $2.7$ s$^{}$\\
$10^{34}$ & $2{,}312{,}042$ & $35$ s$^{}$\\
$10^{35}$ & $25{,}035{,}001$ & $304$ s$^{}$\\
$10^{37}$ & $1{,}544{,}940{,}705$ & $6{,}147$ s$^{}$\\
$10^{38}$ & $9{,}854{,}418{,}621$ & $36{,}520$ s$^{}$\\
\bottomrule
\end{tabular}
\caption{Bounded search above $\omega=10$.}
\end{table}

\section{Verification, and what the results rest on}\label{sec:verif}

The following checks are verified during the computation to ensure a wrong rule never eliminates a node that could have contained a valid friend.
\begin{itemize}
\item \emph{Bound tests} for the exponent bound of Proposition~\ref{prop:cor6}.
  Every integer $N$ satisfies $\sigma(N)/N=\tau$ for its own $\tau=\sigma(N)/N$, and the bound and
  its proof hold for such a target, with the primes of the numerator of $\tau$ in place of $3$.
  So for an arbitrary $N$ with even exponents, a prime $r\mid N$, an arbitrary split of the other
  primes of $N$ into known (fixed, free, or free with an upper bound) and unknown ones, and
  $L_i=U_i$, the program must return $A\ge a_r$. This holds in all of about $58{,}000$ random
  instances: $28{,}192$ with several open terms, $28{,}236$ with a single open term, and $1{,}170$
  with unknowns of very different Wieferich levels, which test the separate bounds $L_i$. In more
  than a quarter of them the bound is attained ($A=a_r$), so the test is not trivially satisfied.
  Of twelve deliberately wrong variants of the bound, nine are detected by this test; the other
  three each remove one unknown prime from the count in Step~2 of the proof, an error that the
  slack of the unknown's own term hides, and they are not detected by it.
\item \emph{Planted targets.} For a number $N_0$ of the right shape, the program is run with the
  target $\I(N_0)$ in place of $9/5$ and a target-independent subset of the rules; it must find
  $N_0$ exactly once. This is done for targets with three and four unknown primes in and above the
  band, with free exponents, at nodes carrying a saturation constraint, and for the configurations
  of Sections~\ref{sec:b7}--\ref{sec:f5}.
\item \emph{Lemma checks} against brute force or an independent computation ($43$ checks for
  Block~8, $45$ for the earlier blocks).
\item \emph{Mutation tests.} $21$ deliberately unsound edits of Block~8 and of the corrected
  closure are all detected; $7$ of them only by lemma checks, since no planted target reaches the
  tail of a region step (it would need an exponent beyond the bound $A$ of
  Rule~\ref{rule:region}, hence an unfactorable target).
\item \emph{Adversarial review.} Five separate reviews were carried out by fresh instances of an AI
  system, each given the code and this paper and asked to find a configuration that
  a rule discards although it contains a solution: Block~7, Block~8, the closure and cofactor
  rules with the correction of Section~\ref{sec:f5}, the bounded search of
  Section~\ref{sec:size}, and arithmetic exactness together with the handling of tasks,
  checkpoints and errors. They found no soundness error. They found a number of places where the
  text did not describe the program exactly, all corrected in this version, and they added
  differential tests: the saturated listing, the deepening and the lower lift against reference
  implementations written from the statements here (several thousand random cases), planted
  targets on the boundaries of the constraints, brute-force checks of step~1 of the partition and
  of the budget cases on $729$ and $472$ random configurations, and $1{,}560$ comparisons of the
  bounded search with an independent brute force. That brute force alone also settles $k=13$ and
  $k=14$ for $B=10^{38}$.
\end{itemize}

\section*{Code availability}\label{sec:code}
The complete PARI/GP source of the certifier and of the bounded search, the test suites (planted
targets, lemma checks and mutation tests), the notes containing the proofs of Blocks~7 and~8, and
the ECPP certificates of Section~\ref{sec:primes} are available as a single archive at
\begin{center}
\url{http://elliptic.group/files/solitary10-gp.tar.bz2}
\end{center}
with SHA-256 checksum
\begin{center}
\small\texttt{3b6e24d0ad87ad304a9bec8afad30480143d536dd27b2bb2220329d1cdf30701}
\end{center}
The \texttt{README.md} file in the archive gives the exact commands for every run reported in
Sections~\ref{sec:comp}--\ref{sec:size}.

\section*{Use of an AI system}
The PARI/GP port, Blocks~7 and~8, their proofs, the tests and this paper were produced with the help of
an AI system  in an extended collaboration with the author.

\section*{Acknowledgements}
This work builds directly on the \texttt{solitary10} project of Caio Comitre Rossi and on the
results of Ward and Thackeray.
We thank H.~R.~Thackeray for correspondence about this work and in particular for confirming
the correction to the hypothesis of Corollary~6 of~\cite{thackeray} discussed after
Proposition~\ref{prop:cor6}.

\end{document}